%% file: j-topologies.tex
\documentclass[a4paper]{amsart}
\usepackage{microtype}%if unwanted, comment out or use option "draft"
\usepackage{calc}
\usepackage{enumerate}
\usepackage{ifthen}
\usepackage{epsfig}
\usepackage{graphicx}
\usepackage{latexsym}
\usepackage{amssymb,amsmath}

\begin{document}
%%%%%%%%%%%%%%%%%%%%%%%%%%%%%%%%%%%%%%%%%%%%%%%%
\title{Counting Finite Topologies}
%------------------------------------------
%\author{E. Fischer}
%\address{Faculty of Computer Science
%Technion-Israel Institute of Technology, Haifa Israel}
%\email{eldar@cs.technion.ac.il}
\author{E. Fischer}
\address{Faculty of Computer Science
Technion-Israel Institute of Technology, Haifa Israel}
\email{eldar@cs.technion.ac.il}
%------------------------------------------
\author{J.A. Makowsky}
\address{Faculty of Computer Science
Technion-Israel Institute of Technology, Haifa Israel}
\email{janos@cs.technion.ac.il}
%------------------------------------------
%%%%amsart
%\author{J.A. Makowsky}
%\address{Faculty of Computer Science
%Technion-Israel Institute of Technology, Haifa Israel}
%\email{janos@cs.technion.ac.il}
%%%%\subjclass{tbc}
%%%% mandatory: Please choose ACM 1998 classifications from http://www.acm.org/about/class/ccs98-html . 
%%%%E.g., cite as "F.1.1 Models of Computation". 
\keywords{Finite topologies, MC-finiteness, congruences, Specker-Blatter Theorem}
%%%% mandatory: Please provide 1-5 keywords
%%%%%%%%%%%%%%%%%%%%%%%%%%%%%%%%%%%%%%%%%%%%%%%%%%%%%%%%%

%Editor-only macros (do not touch as author)
%%%%%%%%%%%%%%%%%%%%%%%%%%%%%%%%%%%
%\serieslogo{}%please provide filename (without suffix)
%\volumeinfo%(easychair interface)
  %{Billy Editor, Bill Editors}% editors
  %{2}% number of editors: 1, 2, ....
  %{TTL-2015, Tools for Teaching Logic }% event
  %{1}% volume
  %{1}% issue
  %{1}% starting page number
%\EventShortName{}
%\DOI{10.4230/LIPIcs.xxx.yyy.p}% to be completed by the volume editor
%%%%%%%%%%%%%%%%%%%%%%%%%%%%%%%%%%%%%%%%%%%%%%%%%%%%%%%%%
%--------------- end lipic----------------------------
%\usepackage{theorem}
%-------------------------

\newcommand{\angl}[1]{\left\langle #1 \right\rangle}
\newcommand{\card}[3]{card_{\mathcal{#1},\bar{#2}}(#3(\bar{#2}))}
\newif\ifmargin
%\margintrue
\marginfalse
\newif\ifshort
\shorttrue
\newif\ifskip
\skiptrue
\newcommand{\NN}{\mathbb{N}}
\newcommand{\ZZ}{\mathbb{Z}}
\newcommand{\FP}{\mathrm{FPT}}
%-----------------------------------------------------------------------------
%\input{e-macros}
\newtheorem{theorem}{Theorem}
\newtheorem{proposition}[theorem]{\bf Proposition}
\newtheorem{examples}[theorem]{\bf Examples}
\newtheorem{example}[theorem]{\bf Example}
\newtheorem{problem}{\bf Problem}
\newtheorem{conclusion}{\bf Conclusion}
\newtheorem{remark}[theorem]{\bf Remark}
\newtheorem{remarks}[theorem]{\bf Remarks}
\newtheorem{definition}{Definition}
\newtheorem{corollary}[theorem]{Corollary}
%---------------------------------------------------
\newtheorem{lesson}{Lesson}
\newtheorem{defi}{Definition}[section]
\newtheorem{conjecture}{Conjecture}
\newtheorem{ex}{Example}[section]
\newtheorem{lemma}[theorem]{Lemma}
\newtheorem{coro}[theorem]{Corollary}
\newtheorem{conj}[theorem]{Conjecture}
\newtheorem{cons}[theorem]{Consequence}
\newtheorem{obs}[theorem]{Observation}
\newtheorem{claim}[theorem]{Claim}
\newtheorem{fact}[theorem]{Fact}
\newtheorem{oproblem}[theorem]{Problem}
% CUSTOMIZING NUMBERED LISTS
\newenvironment{renumerate}{\begin{enumerate}}{\end{enumerate}}
\renewcommand{\theenumi}{\roman{enumi}}
\renewcommand{\labelenumi}{(\roman{enumi})}
\renewcommand{\labelenumii}{(\roman{enumi}.\alph{enumii})}
%----------------------------------------------------------
\renewcommand{\tilde}{\widetilde}
\renewcommand{\bar}{\overline}
%MATH MODE ABBREVIATIONS
\newcommand{\dd}{\mathrm{D}}
\newcommand{\WFF}{\mathrm{WFF}}
\newcommand{\SOL}{\mathrm{SOL}}
\newcommand{\FOL}{\mathrm{FOL}}
\newcommand{\MSOL}{\mathrm{MSOL}}
\newcommand{\CMSOL}{\mathrm{CMSOL}}
\newcommand{\CFOL}{\mathrm{CFOL}}
\newcommand{\TMSOL}{\mathrm{TMSOL}}
\newcommand{\TCMSOL}{\mathrm{TCMSOL}}
\newcommand{\TFOL}{\mbox{TFOL}}
\newcommand{\IFPL}{\mathrm{IFPL}}
\newcommand{\FPL}{\mathrm{FPL}}
\newcommand{\FPT}{\mathrm{FPT}}
\newcommand{\SEN}{\mbox{\bf SEN}}
\newcommand{\WFTF}{\mbox{\bf WFTF}}
\newcommand{\TFOF}{\mbox{\bf TFOF}}
\newcommand{\FOF}{\mbox{\bf FOF}}
\newcommand{\NNF}{\mbox{\bf NNF}}
\newcommand{\N}{{\mathbb N}}
\newcommand{\bN}{{\mathbb N}}
\newcommand{\bR}{{\mathbb R}}
\newcommand{\HF}{\mbox{\bf HF}}
\newcommand{\CNF}{\mbox{\bf CNF}}
\newcommand{\PNF}{\mbox{\bf PNF}}
\newcommand{\QF}{\mbox{\bf QF}}
\newcommand{\DNF}{\mbox{\bf DNF}}
\newcommand{\DISJ}{\mbox{\bf DISJ}}
\newcommand{\CONJ}{\mbox{\bf CONJ}}
\newcommand{\Ass}{\mbox{Ass}}
\newcommand{\Var}{\mbox{Var}}
\newcommand{\Support}{\mbox{Support}}
\newcommand{\V}{\mbox{\bf Var}}
\newcommand{\fA}{{\mathfrak A}}
\newcommand{\fB}{{\mathfrak B}}
\newcommand{\fN}{{\mathfrak N}}
\newcommand{\fZ}{{\mathfrak Z}}
\newcommand{\fQ}{{\mathfrak Q}}
\newcommand{\Aa}{{\mathfrak A}}
\newcommand{\Bb}{{\mathfrak B}}
\newcommand{\Cc}{{\mathfrak C}}
\newcommand{\Gg}{{\mathfrak G}}
\newcommand{\Ww}{{\mathfrak W}}
\newcommand{\Rr}{{\mathfrak R}}
\newcommand{\Nn}{{\mathfrak N}}
\newcommand{\Zz}{{\mathfrak Z}}
\newcommand{\Qq}{{\mathfrak Q}}
\newcommand{\F}{{\mathbf F}}
\newcommand{\T}{{\mathbf T}}
\newcommand{\Z}{{\mathbb Z}}
\newcommand{\R}{{\mathbb R}}
\newcommand{\C}{{\mathbb C}}
\newcommand{\Q}{{\mathbb Q}}
\newcommand{\bP}{{\mathbf P}}
\newcommand{\bPH}{{\mathbf{PH}}}
\newcommand{\bNP}{{\mathbf{NP}}}
\newcommand{\bFPT}{{\mathbf{FPT}}}

\newcommand{\MT}{\mbox{MT}}
\newcommand{\TT}{\mbox{TT}}
\newcommand{\cL}{\mathcal{L}}
\newcommand{\cT}{\mathcal{T}}
\newcommand{\cU}{\mathcal{U}}
\newcommand{\cQ}{\mathcal{Q}}
\newcommand{\cB}{\mathcal{B}}
%-----------------------------------------------------------------------------
\begin{abstract}
\input{j-abstract}
\end{abstract}
\maketitle
\begin{center}
\today
\end{center}
%-----------------------------------------------------------------
\sloppy
%---------------------------------------------
\tableofcontents
\newpage
\input{j-files}

%-----------------------------------------------------------------
\end{document}

%% file: j-abstract.tex
A finite topology $\mathcal{T} = (A, T)$ consists of a finite set $A$ together
with a family $T$ of subsets of $A$, the open sets, satisfying the axioms of a topology. $\mathcal{T}_{\phi}(n)$ is the number of distinct topolgies $T$ of subsets
of $[n]$ which satisfy $\phi$, where $\phi$ is a property of topologies expressible
in $\TMSOL$, topological monadic second order logic.

A sequence $s(n)$ of integers is C-finite if it satisfies a linear recurrence
relation with constant coefficients. It is MC-finite of for every modulus $m$
the sequence $s^m(n) = s(n) \pmod{m}$ satisfies a linear recurrence
relation with constant coefficients depending in $m$. In general $\mathcal{T}_{\phi}(n)$ is not C-finite, because it grows too fast.

In this paper we show that $\mathcal{T}_{\phi}(n)$ is MC-finite for every
$\phi \in \TMSOL$. We also show that this is still true for $\TCMSOL$, the extensions of
$\TMSOL$ with modular counting quantifiers.

%% file: j-files.tex
%\newpage
\input{j-intro}

\input{j-background}
\input{j-logic}

\input{j-counting}
\input{j-proofs}
\input{j-conclu}

%\newpage
\bibliographystyle{plain}
\bibliography{j-ref}
\appendix
\input{j-normal}

%% file: j-intro.tex
\section{Introduction}
In the last decade finite topologies have received renewed attention
due to their role in image analysis and data science.
There is a vast literature testifying to this.
We just mention two references as typical examples, \cite{kovalevsky1989finite,chazal2021introduction}.
For mathematical applications of finite topologies, see \cite{barmak2011algebraic}.
Finite metric spaces are studied in \cite{bourgain1986type,linial2003finite}.
The model theory of topological spaces was studied in the late 1970ies, see 
\cite{makowskytopological,flum2006topological,makowsky1981topological,bk:BFxv}.

In \cite{broder1984r} A. Broder introduced the restricted r-Stirling numbers and r-Bell numbers.
They have found various applications in enumerative combinatorics, e.g. see \cite{benyi2019restricted}.
Inspired by this we study in this paper the number of finite topologies on a finite set subject
to various restrictions.

Assume you are given the set $[r+n] = \{1, 2, 3, \ldots , r+m\}$ and you want to count the number of
finite topolgies on the set $[r+n]$ such that
\begin{enumerate}[(i)]
\item Each of the elements in $[r]$ are in a different component, and
\item Each connected component has odd size.
\end{enumerate}
Or think of variations thereof,
where the topologies on the elements of $[r]$ satisfy some prescribed topological configuration, such as
all the singletons in $[r]$ being closed sets, or being pairwise separable by an open set.
In the applications from \cite{kovalevsky1989finite,chazal2021introduction} counting finite topologies
with such or similar restrictions might be interesting.

Counting finite topologies is a difficult problem.
Even for the case without restrictions no explicit formula is known.
There are some asymptotic results, but the best results known so far are congruences modulo a fixed integer $m$.
Computing $T(n)$ for $n=1,2,3,4$ can be done by hand.
There are two papers giving values for $T(5)$: In \cite{EHL} it says that $T(5)= 6942$.
In \cite{Shaafat} it says that $T(5)= 7181$.
E. Specker in 1980 got interested in modular counting of finite topologies
in order to prove which one of the these claims must be false.

A sequence $s(n)$ of integers is {\em C-finite}, 
if the sequence satisfies a linear recurrence relation with constant coefficients.
A sequence $s(n)$ of integers is {\em MC-finite, modularily C-finite}, if for every integer $m$ the sequence
$s^m(n) =s(n) \pmod{m}$ is an ultimately periodic sequence of positive integers. 
E. Specker showed in \cite{ar:Specker88,specker2011application} the following theorem:
\begin{theorem}
\begin{enumerate}[(i)]
\item
The number of finite topologies $T(n)$ on $[n]$ is not C-finite, but it is MC-finite.
\item
For every $m \in \N^+$ there is a polynomial time algorithm which computes $T(n)$ modulo $m$.
\item
$T(n) = 2 \mod{5}$, hence $T(5) \neq 7181$.
\end{enumerate}
\end{theorem}
The proof of his theorem uses both logic and advanced combinatorics.

The purpose of this paper is to show similar results for the number of topologies with restrictions as suggested above.
In the presence of the restrictions we have in mind, Specker's method cannot be applied directly, 
but it can be applied using our recent results together with a suitable definition of a logic
from \cite{fischer2022extensions}.

We will use logic to make the framework of the restrictions precise. 
Let $\mathcal{T}=(X, \mathcal{U})$ be a finite topological
space on the finite set $X$ and $\cU$ the family of open sets in $X$.
We associate with $\cT$ a two sorted first order structure $\cT' = (X, \cU, E)$
where $E \subseteq X \times \cU$ and $x E U$ says that $x$ is an element of $U$. 
If $X = [r+n]$ we use constant symbols $a_1, \ldots , a_r$ which have a fixed interpretation:
$a_i$ is interpreted by $i \in [r]$. We say that the constant symbols $a_i$ are {\em hard-wired}.
The topological restrictions are now described by first order formulas $\phi(a_1, \ldots, a_r)$
over the structure $([r+n], \cU, \in, a_1, \ldots, a_r)$.
We denote by $T_{\phi, r}(n)$ the number of topologies on the set $[r+n]$ which satisfy $\phi(a_1, \ldots, a_r)$.
For a positive integer $m$
we denote by $T_{\phi, r}^m(n)$ the  sequence $T_{\phi, r}(n)$ modulo $m$.

\subsection{Main result} \ \\
Our main result is stated here for topological first order logic $\TFOL$:
\begin{theorem}
\label{th:main}
\begin{enumerate}[(i)]
\item
For every formula $\phi$ of $\TFOL$ and every positive integer $m$, the sequence
$T_{\phi, r}^m(n)$ is ultimately periodic modulo $m$. In other words $T_{\phi, r}^m(n)$ is MC-finite.
\item
Given $\phi$ and $m$, the sequence $T_{\phi, r}^m(n)$ is Fixed Parameter Tractable ($\bFPT$) where the parameters depend on 
$\phi, r$ and $m$.
\end{enumerate}
\end{theorem}

The proof uses recent results on extensions of Specker's method due to the authors,
\cite{fischer2022extensions}. It also uses model theoretic methods as described in \cite{bk:BFxv}.
The same method was applied to prove congruences for restricted Bell and Stirling numbers in 
\cite{FFMR}.
%\cite{fischer2023mc}.
%\cite{FischerMaRa2023}.

One of our main contributions lies in identifying the logic $\TCMSOL$, 
a topological version of Monadic Second Order Logic with modular counting.
This allows us to prove 
Theorem \ref{th:main-1} in Section \ref{se:proofs}, which is like Theorem \ref{th:main} but stated for $\TCMSOL$
instead of $\TFOL$.

%% file: j-background.tex
\section{Background}
\subsection{C-finite and MC-finite sequences of integers}
\label{se:mcfinite}
A sequence of integers $s(n)$ is {\em C-finite}\footnote{
These are also called constant-recursive sequences
or  linear-recursive sequences in the literature.
}
if there are constants $p, q \in \N$ and $c_i \in \Z, 0 \leq i \leq p-1$ such that for all $n \geq q$ the
linear recurrence relation
$$
s(n+p) = \sum_{i=0}^{p-1} c_i s(n+i)
$$
holds for $s(n)$.
C-finite sequences have limited growth, see e.g. \cite{everest2003recurrence,kauers2011concrete}:
\begin{proposition}
\label{prop:c-finite}
Let $s_n$ be a C-finite sequence of integers. Then there exists $c \in \N^+$ such that for all $n \in \N$,
$a_n \leq 2^{cn}$.
\end{proposition}
Actually, a lot more can be said, see \cite{flajolet2009analytic}, but we do not need it for our purposes.

To prove that a sequence $s(n)$ of integers is not C-finite, we can use Proposition \ref{prop:c-finite}.
To prove that a sequence $s(n)$ of integers is C-finite, there are several methods:
One can try to find an explicit recurrence relation, one can exhibit a rational generating function,
or one can use a method based on model theory as described in 
\cite{fischer2008linear,fischer2011application}.
%The last method will be briefly discussed in
%Section \ref{se:fm} and further explained in Appendix \ref{se:c-finite}. It is referred to as method FM.

A sequence of integers $s(n)$ is {modular C-finite}, abbreviated as {\em MC-finite}, 
if for every $m \in \N$ there are constants $p_m, q_m \in \N^+$ such that
for every  $n \geq q_m$ there is a linear recurrence relation
$$
s(n+p_m) \equiv \sum_{i=0}^{p_m-1} c_{i,m}  s(n+i) \bmod{m}
$$
with constant coefficients $c_{i,m} \in \Z$. 
%and $n \geq n_0$.
Note that the coefficients $c_{i,m}$ and both $p_m$ and $q_m$  generally do depend on $m$.

We denote by $s^m(n)$ the sequence $s(n) \bmod{m}$. 
\begin{proposition}
The sequence $s(n)$ is MC-finite iff $s^m(n)$ is ultimately periodic for every $m$.
\end{proposition}
\begin{proof}
MC-finiteness clearly implies periodicity. The converse is from \cite{reeds1985shift}.
\end{proof}

Clearly, if a sequence $s(n)$ is C-finite it is also MC-finite with $r_m=r$ and $c_{i,m}=c_i$ for all $m$.
The converse is not true, there are uncountably many MC-finite sequences, but only
countably many C-finite sequences with integer coefficients, see Proposition \ref{pr:many} below.
%Here are some typical examples:
\begin{examples}\ 
\label{ex:mc}
\begin{enumerate}[(i)]
\item
The Fibonacci sequence is C-finite.
\item
If $s(n)$ is C-finite it has at most simple exponential growth, by Proposition \ref{prop:c-finite}.
\item
The Bell numbers $B(n)$ are {\em not C-finite}, but are {\em MC-finite}.
\item
Let $f(n)$ be any integer sequence. The sequence $s_1(n)=2\cdot f(n)$ is ultimately periodic modulo $2$,
but not necessarily MC-finite.
\item
Let $g(n)$ be any integer sequence.
%grow arbitrarily fast. 
The sequence  $s_2(n) = n!\cdot g(n)$ is MC-finite.
\label{many-mc}
%We conclude that there are uncountably many monotonously increasing sequences which are MC-finite.
\item
The sequence $s_3(n)= \frac{1}{2} {2n \choose n}$ is not MC-finite: 
$s_3(n)$ is odd iff $n$ is a power of $2$, and otherwise it is even (Lucas, 1878).
A proof may be found in \cite[Exercise 5.61]{graham1989concrete} or in \cite{specker1990application}.
\item
The Catalan numbers $C(n) = \frac{1}{n+1}{2n \choose n}$ are not MC-finite,
since $C(n)$ is odd iff $n$ is a Mersenne number, i.e.,  $n = 2^m-1$ for some $m$,
see \cite[Chapter 13]{koshy2008catalan}.
%For a recent short proof of this, see \cite{KoshySalmassi}, and for an equivalent result from 1973, 
%see \cite[Theorem2]{alter1973binary}.
\ifskip\else
\item
\label{many-nonmc}
Let $p$ be a prime and $f(n)$ be monotone increasing.
The sequence 
$$
s(n) = \begin{cases}
p^{f(n)} & n \neq p^{f(n)} \\
p^{f(n)}+1 & n = p^{f(n)}
\end{cases}
$$
is monotone increasing but  not ultimately periodic modulo $p$, hence not MC-finite.
%We conclude that there are uncountably many monotonously increasing sequences which are not MC-finite.
\fi %skip
\item
\label{many-nonmc}
Let $p$ be a prime and $f(n)$ be monotone increasing.
The sequence $s(n)=p\cdot f(n)+z(n)$, where $z(n)$ is defined to equal $1$ if $n$ is a power of $p$ 
and to equal $0$ for any other $n$,
is monotone increasing but  not ultimately periodic modulo $p$, hence not MC-finite.
%We conclude that there are uncountably many monotonously increasing sequences which are not MC-finite.
\end{enumerate}
\end{examples}

\begin{proposition}
\label{pr:many}
\begin{enumerate}[(i)]
\item
There are uncountably many monotone increasing sequences which are MC-finite, and uncountably many
which are not MC-finite.
\item
Almost all bounded integer sequences are not MC-finite.
\end{enumerate}
\end{proposition}
\begin{proof}
(i) follows from
Examples \ref{ex:mc}
(\ref{many-mc}) and (\ref{many-nonmc}).
(ii) is shown in
Proposition \ref{pr:normal} in  Appendix \ref{se:normal}
\end{proof}

Although we are mostly interested in MC-finite sequences $s(n)$, it would be natural to check in each example
whether the sequence $s(n)$ is also C-finite. In most concrete examples the answer is negative, 
which can be seen by a growth argument. Proposition \ref{pr:C-fin} in Section \ref{se:models} gives a model theoretic
tool for finite topological structures.
However, we will not elaborate this further.  

\subsection{Counting finite topologies}
\label{se:ftop}
Here we follow the presentation from \cite{may2003finite}.
Let $\mathcal{T}=(X, \mathcal{U})$ be a finite topological
space on the finite labeled set $X$ and let $\cU$ be the family of open sets in $X$.
Counting topologies on $X$ is defined as counting the number of distinct families $\cU$ of subsets of $X$.
By Alexandroff's Theorem \ref{th:0}(i)  and \ref{th:1931} 
this is equivalent to counting the number of labeled finite quasi-orders.
Let $T(n)$ and $T_0(n)$ be the number of topologies  and $T_0$-topologies respectively on
a the set $[n] = \{1, \ldots, n \}$.
Recall that a topology on $[n]$ is $T_0$
if for all $a, b \in [n]$, there is some open set
containing one but not both of them. 
No explicit formulas for $T(n)$ and $T_0(n)$ are known. 

The following is known:
\begin{theorem}
\label{th:0}
\ 
\begin{enumerate}[(i)]
\item
$T(n) = Q(n)$, where $Q(n)$ is the number of pre-orders on $[n]$, \cite{may2003finite}.
It is $A000798$ in the Online Encyclopedia of Integer Sequences,
\cite{oeis}.
\item
$T_0(n) = P(n)$, where $P(n)$ is the number of partial orders on $[n]$, \cite{may2003finite}.
It is $A001035$ in the Online Encyclopedia of Integer Sequences.
\item
$Q(n) = \sum_{k=0}^n S(n, k)\cdot P(k)$, where $S(n,k)$ is the Stirling number of the second kind, \cite{butler1973enumeration}.
\item
$B(n) \leq P(n) \leq Q(n)$, where $B(n)$ are the Bell numbers, which count the number of equivalence relations
on $[n]$. Furthermore, see \cite{de1981asymptotic,berend2010improved},
$$
\left(\frac{n}{e \ln n}\right)^n \leq B(n) \leq \left(\frac{n}{e^{1-\epsilon} \ln n}\right)^n,
$$
\item
The logarithm with base $2$ of both $T(n)$ and $T_0(n)$
goes asymptotically to $\frac{n^2}{4}$
as $n$ goes to infinity,
\cite{kleitman1970number}. 
\end{enumerate}
\end{theorem}

\begin{theorem}
\begin{enumerate}[(i)]
\item
$T(n)$ and $T_0(n)$ are not C-finite.
\item
$T(n)$ and $T_0(n)$ are MC-finite.
\end{enumerate}
\end{theorem}
\begin{proof}
(i) follows from Theorem \ref{th:0}(iv).
\\
(ii) follows from Theorem \ref{th:0}(i) and (ii) and the Specker-Blatter Theorem \ref{th:SB}.
\end{proof}

%% file: j-logic.tex
\section{Topologies as relational structures}
\label{se:models}
\subsection{The logics $\TCMSOL$ and $\CMSOL$}

Let $\mathcal{T}=(X, \mathcal{U}, \in )$ be a finite topological space 
on the finite set $X$ and let $\mathcal{U}$ be the family of open sets in $X$. The relation $\in$ is the
natural membership relation between elements of $X$ and sets in $\mathcal{U}$.
We associate with a finite topological space $(X, \mathcal{U})$ a two sorted relational structure $\cT = (X, \cU, E)$,
where $E \subseteq X \times \cU$ is a binary relation and $E(x,U)$ says that $x$ is an element of $U$. 
$E$ is required to satisfy the extensionality axiom 
$$
\forall U, V \in \cU 
(u=v \leftrightarrow 
\forall x \in X (E(x,U) \leftrightarrow E(x,v))
).
$$
There is a natural bijection between finite topological spaces and and their associated first order structures.
Given $\mathcal{T}=(X, \mathcal{U}, \in )$, a finite topological space, we define the first order structure
$t((X, \cU, E ))$
by setting $\cU = \mathcal{U}$ and $E(x,U)$ iff $x \in U$.
Conversely, given $\cT = (X, \cU, E)$ which satisfies the extensionality axiom, we 
define $t^{-1}((X, \cU, E))$ by setting for $U \in  \mathcal{U}$
$$
x \in U \text{  iff  } E(x,U).
$$
\begin{proposition}
\label{pr:t-bijection}
$t$ is a bijection between finite topological spaces and and their associated first order structures.
Furthrmore, $t^{-1}$ is its inverse.
\end{proposition}

We denote by $\TCMSOL$ the monadic second order logic for structures of this form possibly augmented by constant symbols.
We allow quantification over {\em subsets of $X$}, but only quantification over {\em elements of $\; \cU$}.
Furthermore, we have a modular counting quantifier $C_{m,a}x \phi(x)$ which says that modulo $m$ there are $a$
elements satisfying $\phi(x)$. $\TFOL$ is the logic without second order quantification and without modular
counting. Similarly, $\CMSOL$ is defined as $\TCMSOL$ for one-sorted structures with one binary relation.

For {\em finite structures} of the form $\cT$ the following are $\TCMSOL$-definable.
\begin{enumerate}[(i)]
\item
$\cU$ is a topology for the finite set $A$: 
(i) $\emptyset  \in \cU$, $A \in \cU$.
(ii) $\cU$ is closed under unions.
(iii) $\cU$ is closed under intersections.
\item
$\cU$ is $T_0$: $\forall a,b \in A \exists U \in \cU 
(E(a,U) \wedge \neg E(b,U) 
\vee
(\neg E(a,U) \wedge  E(b,U)) 
$.
\item
$\cU$ is $T_1$: $\forall a \in A (A - {a} \in \cU)$.
\item
$X$ is connected: There are no two non-empty disjoint open sets $U_1, U_2$ with $U_1 \cup U_2 =X$.
\item
The $\TFOL$-formula  $\phi_{U_x}(x,U)$  says that $U$ is the smallest open set containing $x$: 
\\
%$U_x \in \cU \wedge E(x,U_x) \wedge (\forall U \in \cU (E(x,U) \rightarrow U_x \subseteq U)$.
%---
$(U\in\mathcal{U})\wedge (E(x,U))\wedge (\forall V\in\mathcal{U}(E(x,V)\to V\subseteq U))$
%----
\item
A typical formula which is in $\TCMSOL$ would be: There is a set of points of even cardinality
which is not an open set.
\end{enumerate}

\subsection{Hard-wired constant symbols}
Let $\bar{a}=(a_1, \ldots, a_k)$ be $k$ constant symbols.
For each of them there are $n$ possible interpretations in the set $[n]$.
However, we say that $(a_1, \ldots, a_r)$, for $r \leq k$ are
{\em hard-wired} on $[n]$, if $a_i$ is interpreted by $i \in [n]$.
In the presence of constant symbols (hard-wired or not) $a_1, \ldots, a_k$ we can say:
\begin{enumerate}[(i)]
\item
%$\forall U \bigwedge_i^r E(a_i, U)$, i.e., they form a minimal non-empty open set.
%---
$\{a_1, \ldots, a_k\}$ is a minimal non-empty open set:
\begin{gather}
\exists U\in\mathcal{U}((\bigwedge_{i=1}^r E(a_i,U))
\notag \\
\wedge 
%\notag \\
(\forall x(\bigwedge_{i=1}^r (x\neq a_i)\to\neg E(x,U)))
\notag \\
\wedge
%\notag \\
(\forall V\in\mathcal{U}((V\subseteq U)\to((V=U)\vee(V=\emptyset)))))
\notag
\end{gather}
%(maybe it can be "abbreviated" a bit, about half of it just means $U=\{a_1,\ldots,a_r\}$)
%---
\item
There are pairwise disjoint open sets $ U_1, \ldots , U_r$ such that $a_i$ in $U_i$.
\item
The elements denoted by $a_i$ are all in different connected components.
\end{enumerate}
In analogy to Broder's $r$-Stirling numbers, we also count finite topologies restricted by $\TCMSOL$-formulas
with hard-wired constant symbols.

\subsection{Topological model theory}
In his retiring presidential address presented at the annual meeting of the Association
for Symbolic Logic in Dallas, January 1973
Abraham Robinson suggested, among other topics, to develop a model theory for topological structures,
\cite{robinson1973metamathematical}.
This led to several approaches described in
\cite{ziegler1976language,flum2006topological,bk:BFxv,makowsky1981topological,makowskytopological}.
M. Ziegler introduced the logic $L_t$, which is a fragment of $\TFOL$ with the additional property that
it is {\em basis-invariant} in the following sense:

Let
$\phi$ be a formula of $L_t$ and $\mathcal{T}=(X, \mathcal{U})$ be a (not necessarily finite) 
topological structure, and let $\mathcal{B}$ be any basis for $\mathcal{U}$.
%$\mathcal{T}=(X, \mathcal{U})$.
Then 
$$\mathcal{T}=(X, \mathcal{U}) \models \phi \text{ iff  } \mathcal{T}=(X, \mathcal{B}) \models \phi.$$
%for every basis $\mathcal{B}$ of $\mathcal{U}$.
In fact in \cite{ziegler1976language,flum2006topological,bk:BFxv} $X$ can be replaced by arbitrary first order
structures. $L_t$ now shares most model theoretic characterstics of first order logic, like compactness, L\"owenheim-Skolem
theorems, preservation theorems, etc.  However, if we restrict the topological structures for $L_t$ to be finite,
this is not true anymore.
A topological structure $\mathcal{T}=(X, \mathcal{U})$ is an {\em open substructure} of
$\mathcal{T}'=(Y, \mathcal{V})$ if $X \subseteq Y$, $\mathcal{U} = \{A \subset X: A = A' \cap X, A' \in \mathcal{V})$,
and $X \in \mathcal{V})$. We write $\mathcal{T} \subseteq_o \mathcal{T'}$.
A formula $\phi \in L_t$ is {\em preserved under open extensions}
if for every  pair of topological structures $\mathcal{T} \subseteq_o \mathcal{T}'$ we have
$$
\mathcal{T} \models \phi \text{  implies  } \mathcal{T}' \models \phi .
$$
In \cite{flum2006topological} there is also a syntactical characterization of the formulas preserved under open extensions.
However, it fails if restricted to finite structures.
Nevertheless, we can use this to show (here without proof):

\begin{proposition}
\label{pr:C-fin}
Let $\phi \in L_t$ with $r$ constant symbols (hard-wired or not) which  has arbitrarily large finite models and 
is preserved under open extensions.
Then $T_{\phi,r}(n)$ is not C-finite.
\end{proposition}

%% file: j-counting.tex
\section{Most general Specker-Blatter Theorem}
\label{se:SB}
%\section{Counting restricted topologies: MC-finiteness}
In order to prove  Theorem \ref{th:main} for $\TCMSOL$,
we now state the Specker-Blatter Theorem for $\CMSOL$ with hard-wired constants,
\cite{fischer2022extensions}. This generalizes substantially the original  Specker-Blatter Theorem from 1981,
%\ref{th:specker}, 
but is
still formulated for one sorted relational structures for binary and unary relation symbols.
Note that in \cite{fischer2022extensions} it is also shown that Theorem \ref{th:SB}
does not hold for ternary relations. 

Let $\tau$ be a finite set of binary and unary relation symbols.
Let $Sp_{\phi, r}(n)$ be the number of labeld $\tau$-structures on the set $[r+n]$
where the elements of $[r]$ are hard-wired and which satisfy the formula $\phi$ of $\CMSOL$.

\begin{theorem}[Most general Specker-Blatter Theorem]
\label{th:SB}
\begin{enumerate}[(i)]
\item
$Sp_{\phi,r}(n)$ is MC-finite.
\item
For every $m \in \N^+$ the sequence $Sp_{\phi,r}^m(n) = Sp_{\phi,r}(n) \bmod{m}$
is computable in polynomial time. In fact it is in $\FPT$ (Fixed Paramater Tractable) 
with parameters $\phi, r$ and  $m$.
\end{enumerate}
\end{theorem}
For the technical details and the history of this theorem, the logically inclined reader should consult
\cite{fischer2022extensions}. 

The theorem holds for a fixed number of hard-wired constants. It does not hold for a hard-wired relation.
A relation $R \subseteq [n]^s$ is hard-wired if $R$ has a unique interpretation on $[n]^s$.
We can view the natural order $NO \subseteq [n]^2$ as hard-wired.
The number of equivalence relations (set partitions) on $[n]$ is given by the Bell numbers $B(n)$ which are MC-finite.
They satisfy the hypothesis of the Specker-Blatter Theorem, as an equivalence relation is definable in $\FOL$.

\begin{proposition}
The Specker-Blatter Theorem does not hold in the presence a hard-wired linear order on $[n]$.
\end{proposition}
\begin{proof}
Let $A$ and $B$ be two blocks of a partition of $[n]$.
$A$ and $B$ are {\em crossing} if there are elements
$a_1, a_2 \in A$ and $b_1, b_2 \in B$ such that 
$a_1 < b_1 < a_2 < b_2$ or
$b_1 < a_1 < b_2 < a_2$.
The number $B(n)^{nc}$ of 
non-crossing set partitions on $[n]$ is one of the interpretations of the Catalan numbers $C(n)$, 
\cite{roman2015introduction}, hence $C(n) =B(n)^{nc}$.
Non-crossing set partitions are definable in $\FOL$ in presence of the hard-wired natural order on $[n]$.
But the Catalan numbers are not MC-finite.
\end{proof}

%Our next task is now to show how to reduce the theorem for $\TCMSOL$ to Theorem \ref{th:SB} above.

%% file: j-proofs.tex
%revision started  August 10, 2023
\section{Proof of Theorem \ref{th:main}}
\label{se:proofs}
\subsection{Alexandroff's Theorem}
Our next task is now to show how to reduce the theorem for $\TCMSOL$ to Theorem \ref{th:SB} of Section \ref{se:SB}.
Besides using Theorem \ref{th:SB} we use several classical facts about finite topologies, taken from \cite{may2003finite}.

First we state a bijection theorem. 

\begin{theorem}[Alexandroff, 1931]
\label{th:1931}
There are 
bijections $\alpha$ and $\alpha'$ 
between finite topologies and finite quasi-orders
and vice versa. Furthermore, they are the inverses of each other.
\end{theorem}

\begin{proof}
Let $\cT =(X, \cU, E)$ be a finite topology. For $x \in X$ define $U_x$ to be intersection of all open sets
which contain $x$. The sets $U_x$ form a basis for $\cU$.
Define a relation $\leq_{\cU}$ on the set $X$ by $x \leq_{\cU} y$
if $x \in U_y$ or, equivalently, $U_x \subseteq U_y$ . Write $x <_{\cU} y$ if the inclusion is proper.
The relation $\leq_{\cU}$ is transitive and reflexive, 
hence it defines for each $\cU$ a unique quasi-order $\leq_{\cU}$.
We define $\alpha(\cT) = (X, \leq_{\cU})$. 

Conversely, let $(X, \leq)$ be a quasi-order. The sets $U_x = \{ y \in X: y \leq x \}$ 
form a basis for a topology $\cU_{\leq}$.
We define $\alpha'((X, \leq)) = (X, \cU_{\leq})$. 
\end{proof}

\subsection{Translation schemes}
The maps $\alpha$ and $\alpha'$ from Theorem \ref{th:1931}  are actually definable in $\MSOL$ and $\TMSOL$ respectively
by translation schemes as defined below.
This is needed to construct an algorithm which translates a formula $\theta$ of $\TCMSOL$ 
into a formula $\theta^{\sharp}$ of $\CMSOL$ over quasi-orders such that 
for every finite quasi-order $\mathfrak{O}$
$$
\alpha'(\mathfrak{Q}) \models \theta \text{ iff } \mathfrak{Q} \models \theta^{\sharp}.
$$
\begin{conclusion}
The number of finite topolologies on a set with $n$ elements which satisfy $\theta$ equals the number of
finite quasi-orders on a set with $n$ elements which satisfy $\theta^{\sharp}$.
\end{conclusion}

We first show the definability of $\alpha$.
Let $\cT =(X, \cU, E)$ be a finite topology.
We want to define inside $\cT$ a quasi-order $\cQ =(X, \leq)$.
For this we exhibit two formulas $\phi(x)$ and $\phi_{\leq}(x,y)$ in $\TMSOL$ which 
form a translation scheme $\Phi =( \phi(x), \phi_{\leq}(x,y))$. 

With $\Phi$ we associate two maps. $\Phi^{\star}$ and $\Phi^{\sharp}$.
$\Phi^{\star}$ maps finite topological spaces into quasi-orders, and
$\Phi^{\sharp}$ maps formulas of $\CMSOL$ into formulas of $\TCMSOL$.
Let $\cT$ a finite topology and $\theta$ a formula of $\CMSOL$ for quasi-orders.
These two maps satisfy the following:
$$
\Phi^{\star}(\cT) \models \theta \text{  iff  }
\cT \models \Phi^{\sharp}(\theta)
$$
and
$$
\alpha(\cT) = \Phi^{\star}(\cT).
$$

For $\alpha'$ we proceed similarily.
Let $\cQ$ be a quasi-order.
We want to define inside $\cQ$ a finite topology using formulas
of $\CMSOL$ given by 
$$
\Psi = (\psi(x), \psi_{open}(U), U(x))
$$
\begin{itemize}
\item
$\psi(x):= x=x$ defines the universe of the topology.
\item
$\psi_{basic}(x, U): = (U(x) \wedge  
(\forall y  (U(y)  \leftrightarrow y \leq x)
)$
defines the basic open sets.
\item
$\psi_{open}(U):= \forall z 
%(U(z) \leftrightarrow   \exists V \exists y V(y) ( \psi_{basic}(y, W) \wedge W(z))$ 
(U(z) \leftrightarrow   \exists V \exists y ( \psi_{basic}(y, V) \wedge V(z))$ 
says that $U$ is a union of basic open sets.
\end{itemize}

With $\Psi$ we associate two maps. $\Psi^{\star}$ and $\Psi^{\sharp}$.
$\Psi^{\star}$ maps finite quasi-orders into finite topologies, and
$\Psi^{\sharp}$ maps formulas of $\TCMSOL$ into formulas of $\CMSOL$.
Let $\cQ$ a finite quasi-order and $\theta$ a formula of $\TCMSOL$.
These two maps satisfy the following:
$$
\Psi^{\star}(\cT) \models \theta \text{  iff  }
\cT \models \Psi^{\sharp}(\theta)
$$
and
$$
\alpha'(\cQ) = \Psi^{\star}(\cQ).
$$

$\Psi^{\sharp}$ is defined inductively.

The atomic formulas of $\TCMSOL$ are
$E(x,U)$, which is translated as $\psi_{open}(U) \wedge U(x)$, and $U(x)$,
which is translated by $U(x)$.
The translations commute with the boolean oeprations and quantifications of $\TCMSOL$.

Let $U_x$ be the intersection of all open sets $U$ which contain $x$.

This can be expressed in $\TFOL$ by the formula
$\phi_{U_x}(x,U)$ from the previous section.
%$$
%\phi_{U_x}(x,U): 
%E(x,U) \leftrightarrow 
%\left[ 
%\forall V (E(x,V) \rightarrow 
%( \forall z E(z, U)  \rightarrow E(z, V))
%\right]
%$$
Now $x \leq y$ can be defined by $U_x \subseteq U_y$, which be expressed as
$$
\phi_{\leq}(x,y):
\forall z  
\left[
\exists U (E(z, U) \wedge \phi_{U_x}(x,U))
\rightarrow
\exists V (E(z, V) \wedge \phi_{U_x}(y,V))
\right]
$$
The translation scheme $\Phi = (x=x, \phi_{\leq}(x,y))$ consists of two formulas.
The first is a tautology in the free variable $x$ and defines the new universe, which in this case is also $X$. 
The second formula, $\phi_{\leq}(x,y))$ defines the quasi-order.

$\Phi$ induces to maps, $\Phi^{\star}$ which maps finite topologies onto quasi-orders over the same universe,
and $\Phi^{\sharp}$, which maps $\CMSOL$-formulas into $\TCMSOL$-formulas, by replacing each occurrence
of $x_1 \leq x_2$ by $\phi_{\leq}(x_1,x_2)$.

\begin{lemma}
$\alpha(\cT) = \Phi^{\star}(\cT)$.
\end{lemma}

In the other direction, let $ \cQ =(X, \leq)$ be a finite quasi-order.
We want to define inside $\cQ$ a topology $\cT =(X, \cU, E)$.
We actually define a structure
$$
\cT' =(X, P(X), \cU, \cB, E, E_{top}, E_{basis})
$$ 
where $P(A)$ is the powerset of $A$, and $\cB$ is a minimal basis
for the topolgy $\cU$.
$E,  E_{top}, E_{basis}$ are the membership relations for elements of $X$ and $P(A), \cU, \cB$ respectively.

Again $X$ can be defined by $x=x$, and $P(A)$ can be defined by $\phi_{set}(X): \forall x (X(x) \leftrightarrow X(x))$,
or, for that matter, by any appropriate tautology in one free variable.

The non-empty basic sets are defined by 
%$$ \cB = \{ B \in P(X): \forall y (y \in B \leftrightarrow \exists x  (y\leq x)\}.$$
$$ \cB = \{ B \in P(X): \exists x \forall y (y \in B \leftrightarrow (y\leq x)\}.$$
Hence we put
%$$\phi_{basis}(B)(y,B): \exists y (B(y) \leftrightarrow \exists x  (y\leq x).$$
%$$\phi_{basis}(B): \forall y (E(y, B) \leftrightarrow \exists x  (y\leq x).$$ 
$$\phi_{basis}(B): \exists x \forall y (E(y, B) \leftrightarrow (y\leq x).$$ 

Then the non-empty open sets are defined by
$$ \cU = \{ U \in P(X): \forall x  (x\in U \leftrightarrow (\exists B \in \cB (x \in B \wedge B \subseteq U))) \}.$$
Hence we put
$$\phi_{top}(U): \forall x  (E(x,U) \leftrightarrow  (\exists B \in \cB (x \in B \wedge B \subseteq U))).$$
%(\exists y \exists B (B(y) \wedge (\forall z (B(z) \rightarrow U(z)))).$$
The translation scheme is now defined by
$$\Psi  = ( x=x, \phi_{set}(X), \phi_{top}(U), \phi_{basis}(B), X(x) ).$$ 

$\Psi$ induces to maps, $\Psi^{\star}$, which maps finite quasi-order onto topologies over the same
underlying set. and $\Psi^{\sharp}$, which transforms $\TCMSOL$-formulas 
%into $\TCMSOL$-formulas,
by replacing each occurrence of $E(x, X), E(x, U), E(x,B)$ by its definitions.

\begin{lemma}
$\alpha'(\cQ) = \Psi(\cQ)$.
\end{lemma}
%§Recall that $\alpha$ is the bijection between finite topologies and quasi-orders from Theorem \ref{th:1931}.
\begin{theorem}
\label{th:trans}
The translation schemes $\Phi$ and $\Psi$ satisfy the following:
\begin{enumerate}[(i)]
\item
$\Phi^{\star}(\cT)= \alpha(\cT)$ and
$\Psi^{\star}(\cQ)= \alpha'(\cQ)$; 
\item
for every 
%$\theta \in \CMSOL^h$ 
$\theta \in \CMSOL$ 
and very finite quasi-order $\cQ$,
$$\alpha(\cT) = \cQ \models \theta  \text{  iff   }  \cT \models \Phi^{\sharp}(\theta).$$
\item
for every 
%$\sigma \in \TCMSOL^h$ 
$\sigma \in \TCMSOL$ 
and very finite topology $\cT$,
$$\alpha'(\cQ) = \cT \models \sigma  \text{  iff   } \cQ \models \Psi^{\sharp}(\sigma).$$
\end{enumerate}
\end{theorem}

%\begin{theorem}[E. Fischer and J.A. Makowsky, 2022]
\begin{theorem}[\cite{fischer2022extensions}]
\label{th:FiMa}
Let $\theta(a_1, \ldots, a_r)$ be a sentence in 
%$\CMSOL^h$
$\CMSOL$
 with $r$ constant symbols, 
and let $S(n)= S_{\theta(a_1, \ldots, a_r)}(n)$ be the number of relations $R \subseteq [n]^2$,
such that 
$$([n], R, (a_1, \ldots, a_r)) \models \theta(a_1, \ldots, a_r).$$ 
%$([n], R, (a_1, \ldots, a_r)) \models \theta(a_1, \ldots, a_r)$. 
In both cases, where the constant symbols are
hard-wired or not, $S$ is MC-finite.
\end{theorem}

Let $\sigma(a_1, \ldots, a_r)$ be a sentence in 
%$\TCMSOL^h$ 
$\TCMSOL$ 
with $r$ constant symbols, 
and let $S^t(n)= S^t_{\sigma(a_1, \ldots, a_r)}(n)$ be the number of topologies on $[n]$
such that 
$$([n], \cU, \, (a_1, \ldots, a_r)) \models \sigma(a_1, \ldots, a_r).$$ 
and let $\hat{S}^t(n)= S^t_{\neg \sigma(a_1, \ldots, a_r)}(n)$ be the number of topologies on $[n]$
such that 
$$([n], \cU, \, (a_1, \ldots, a_r)) \models \neg \sigma(a_1, \ldots, a_r).$$ 
\begin{theorem}
\label{th:main-1}
In both cases, whether the constant symbols are
hard-wired or not, 
\begin{enumerate}[(i)]
\item
$S^t(n)$ 
and
$\hat{S}^t(n)$ 
are MC-finite, but 
\item
at least one of them is not C-finite.
\end{enumerate}
\end{theorem}

%\begin{proof}
%\end{proof}

%% file: j-conclu.tex
\section{Conclusions}
\label{se:conclu}

We have shown that counting finite topologies $T_{\phi, r}(n)$ on a set of $n$ elements 
subject to restrictions with a fixed finited number $r$ of (hard-wired)
constants expressed by $\phi \in \TCMSOL$ is an MC-finite sequence.
The underlying set $X$ of the topology can be equipped with unary and binary relations, 
but counting now also counts the number of their interpretations.
As we have seen in Section \ref{se:SB}
the Specker-Blatter Theorem does not work for hard-wired relations.

No explicit formulas for $T(n)$ and $T_0(n)$ are known, but we have
$$T(n) = \sum_{k=0}^n S(n, k)\cdot T_0(k),$$ 
where $S(n,k)$ is the Stirling number of the second kind. 
%\cite{butler1973enumeration}.

\begin{problem}
Find better descriptions of $T(n)$ and $T_0(n)$.
\end{problem}

The logic $\TCMSOL$ can express many topological properties.  Here are some possibly challenging
test problems.

In analogy to counting various set partitions
on $[n]$ as described in \cite{FFMR},
one can look at topological set partitions of finite topological spaces on $[n]$ 
such that each block is a connected component, an open, or a closed set.
Note that the topology in a topological set partition is not hard-wired.
We count the set partitions of $X$ and the topologies separately. 
Given $[n]$ the Bell numbers $B(n)$ count the partitions of $[n]$ and for each such partition we count the number of topologies
such that the blocks are connected, open or closed.
The Stirling numbers of the second kind $S(n,k)$ count the partitions of $[n]$ into $k$ blocks, and again we can require
that the blocks are connected, open or closed.

Similarlily, we can look at topological spaces on $[r+n]$ such that the hard-wired constants
of $[r]$ are all in different blocks. Like with the Bell and Stirling numbers, 
the number of non-empty blocks may be arbitrary or fixed to $k$ blocks.
If the topology is assumed to be discrete, we get the
resricted Bell numbers $B_r(n)$ and Stirling numbers of the second kind $S_r(n,k)$ from \cite{broder1984r}.
Otherwise, these give different topological versions of $B_r(n)$ and $S_r(n,k)$.

\begin{problem}
What can we say beyond that the number of such topological set partitions
form an MC-finite sequence?
Are some of them C-finite?
\end{problem}

%For every formula $\phi$ of $\TFOL$ and every positive integer $m$, the sequence
%$T_{\phi, r}^m(n)$ is ultimately periodic modulo $m$. In other words $T_{\phi, r}^m(n)$ is MC-finite.

%% file: j-normal.tex
\section{Normal sequences}
\label{se:normal}
Let $s(n)$ be an integer sequence, and $b \in \N^+$.
The sequence $s^b(n) = s(n) \bmod{b}$ 
is normal, if, when we chunk it into substrings of length 
$\ell \ge 1$, then each of the $b^\ell$ possible strings of $[b]^\ell$ appear in $s^b(n)$ with equal limiting frequency. 
It is {\em absolutely normal} if it is normal for every $b$.
The sequence $s^b(n) = s(n) \bmod{b}$ can be viewed as a real number $r_b$ written in base $b$.
A classical theorem  from 1922 by E. Borel says that almost all reals are absolutely normal, 
\cite{everest2003recurrence}.
The proposition below shows that MC-finite integer sequences are very rare.

Let $PR_b$ be the set of integer sequences $s^b(n)$ 
with $s^b(n) = s(n) \bmod{b}$ for some integer sequence $s(n)$.
$PR_b$ is the projection of all integer sequences to sequences over $\Z_b$.
We think of $PR_b$ as a set of reals with the usual topology and its Lebesgue measure.
Let $UP_b \subseteq PR_b$ be the set of sequences $s^b(n) \in PR_b$ which are ultimately periodic.

\begin{proposition}
\label{pr:normal}
\begin{enumerate}[(i)]
\item
Almost all reals are absolutely normal. 
\item
$s(n)$ is MC-finite iff for every $b \in \N^+$ the sequence $s^b(n)$ is ultimately periodic.
\item
If $s^b(n)$ is normal for some $b$, then $s(n)$ is not MC-finite. 
\item
$UP_b \subseteq PR_b$ has measure $0$.
\end{enumerate}
\end{proposition}
Proving that a specific sequence is normal is usually very difficult.
%wikipedia  https://en.wikipedia.org/wiki/Normal_number
It has been an elusive goal to prove the normalcy of numbers that are not artificially constructed. 
While $\sqrt{2}$, $\pi$, $ln(2)$ and $e$ are strongly conjectured to be normal, 
it is still not known whether they are normal or not. 
It has not even been proven that all digits actually occur infinitely many times in the 
decimal expansions of those constants. 
%(for example, in the case of π, the popular claim "every string of numbers eventually occurs in π" is not known to be true).[15] 
It has also been conjectured that every irrational algebraic number is absolutely normal, 
and no counterexamples are known in any base. However, no irrational algebraic number has been proven to be normal in any base. 

\ifskip\else
Here is a challenge:
\begin{conjecture}
The binary sequence $\beta(n) = a(n) \bmod 2$ from Theorem \ref{thm:main} is normal with $b=2$. 
\end{conjecture}
\fi %skip